\documentclass{elsart}

\usepackage{amssymb,amsmath}
\usepackage{graphicx}

\begin{document}

\begin{frontmatter}


% use the thanksref command within \title, \author or \address for
% footnotes;
% use the corauthref command within \author for corresponding author
% footnotes;

\title{Simple formulas for lattice paths avoiding
certain periodic staircase boundaries}


\author[add:uk]{Robin J. Chapman},
\author[add:ccr]{Timothy Y. Chow},
\author[add:ccr]{Amit Khetan},
\author[add:ccr]{David Petrie Moulton},
\author[add:uk]{Robert J. Waters}
\address[add:ccr]{Center for Communications Research, 805 Bunn Drive,
Princeton, NJ 08540, U.S.A.}
\address[add:uk]{Department of Mathematics, University of Bristol,
University Walk, Bristol, BS8~1TW, U.K.}

\begin{abstract}
There is a strikingly simple classical formula for the
number of lattice paths avoiding the line $x=ky$ when
$k$ is a positive integer.
We show that the natural generalization of this simple formula 
continues to hold when the line $x=ky$ is replaced by
certain periodic staircase boundaries---but only under special conditions.
The simple formula fails in general,
and it remains an open question
to what extent our results can be further generalized.
\end{abstract}

\begin{keyword}
ballot sequence \sep zigzag \sep stairstep \sep
touching \sep crossing \sep tennis ball
% \PACS 
\end{keyword}
\end{frontmatter}

\def\fkn{{{\mathcal F}_k(n)}}
\def\gkn{{{\mathcal G}_k(n)}}
% main text
\section{Background and main results}
\label{sec:intro}

Throughout this paper, a {\em lattice path\/} will mean
a lattice path in the plane whose only allowable steps
are north $(0,1)$ and east $(1,0)$.

It is a classical theorem \cite{Bar}\cite{Ber}
that if $k$ is a positive integer,
then the number of lattice paths from
$(0,0)$ to $(a+1,b)$ (where $a\ge kb$) that avoid touching or crossing 
the line $x=ky$ except at $(0,0)$ is given by the formula
\begin{equation}
\label{eq:classic}
\binom{a+b}{b} - k \binom{a+b}{b-1}.
\end{equation}
In fact, more is true: There are
\begin{equation}
\label{eq:classicnw}
\binom{a }{ c-1}\binom{b }{ c-1} - k \binom{a-1}{ c-2}\binom{b+1}{ c}
\end{equation}
such paths with $c-1$ northwest corners.\footnote{By
a northwest corner we mean a corner formed by a north step
followed by an east step.
There is a similar-looking formula 
for paths with a given number of southeast corners.}
This stronger result appears explicitly in~\cite{Kra}
and implicitly even earlier,
but our favorite proofs of all these facts
are the bijective proofs of Goulden and Serrano~\cite{GS}.

It is natural to ask if there are similar simple formulas for
lattice paths from $(0,0)$ to $(a,b)$ that avoid the line $x=ky$,
if $k$ is allowed to be an arbitrary positive rational number.
While one can write down a determinantal formula
(indeed, a determinantal formula exists
for an arbitrarily shaped boundary),
nothing as simple as (\ref{eq:classic}) is known,
and empirical investigation does not suggest any obvious conjecture.

Our first main result is that
for certain periodic staircase boundaries
(instead of straight-line boundaries),
there {\em are\/} simple enumerative formulas
that generalize (\ref{eq:classic}) and (\ref{eq:classicnw}),
at least for certain special starting and ending points.

\begin{defn}
Given positive integers $s$ and $t$,
let $A_{s,t}$ be the infinite staircase path that starts at $(0,t)$,
then takes $s$ steps east, $t$ steps north, $s$ steps east,
$t$ steps north, and so on.
\end{defn}

\begin{defn}
Given a set $S$ of (finite) lattice paths,
take each path $\pi\in S$, and augment it by
prepending a north step to the beginning of~$\pi$
and appending a north step to the end of~$\pi$.
Let $S^+$ denote the resulting set of lattice paths.
\end{defn}

\begin{thm}
\label{thm:corners}
Let $s$, $t$, $n$, and $c$ be positive integers.
\begin{enumerate}
\item Let $S_1$ be the set of lattice paths from
$(0,0)$ to $(sn+1,tn)$ that avoid $A_{s,t}$.  There are
\begin{equation}
\label{eq:nwone}
t \binom{sn}{c-1}\binom{tn}{c-1} - s \binom{sn-1 }{c-2} \binom{tn+1}{c}
\end{equation}
paths in~$S_1^+$ with $c$ northwest corners
(equivalently, $c$ southeast corners).
\item Let $S_2$ be the set of lattice paths from
$(1,0)$ to $(sn,tn-1)$ that avoid $A_{s,t}$.  There are
\begin{equation}
\label{eq:nwtwo}
t \binom{sn-1}{c-1}\binom{tn-1}{c-1} - s \binom{sn-2}{c-2} \binom{tn}{c}
\end{equation}
paths in~$S_2^+$ with $c$ northwest corners
(equivalently, $c$ southeast corners).
\end{enumerate}
\end{thm}

The equivalence between counting northwest and southeast corners
follows because in a lattice path that starts with a north step
and ends with a north step,
the first corner must be a northwest corner and the last corner
must be a southeast corner, and northwest and southeast corners
must alternate.
Also, since $|S| = |S^+|$ for any~$S$,
summing over all $c$ and applying Vandermonde convolution
immediately yields the following corollary.

\begin{cor}
\label{cor:total}
Let $s$, $t$, and $n$ be positive integers.  Then
\begin{equation}
\label{eq:totalone}
|S_1| = t \binom{sn+tn}{tn} - s \binom{sn+tn }{tn-1}
\end{equation}
and
\begin{equation}
\label{eq:totaltwo}
|S_2| = t \binom{sn+tn-2 }{tn-1} - s \binom{sn+tn-2 }{tn-2}.
\end{equation}
\end{cor}

Note that avoiding $A_{k,1}$ is the same as
avoiding $x=ky$ except at $(0,0)$,
so our results generalize (\ref{eq:classic}) and (\ref{eq:classicnw}) 
in one direction, by allowing arbitrary $s$ and~$t$,
but are simultaneously more special in another direction,
since only certain special endpoints are allowed.
More precisely, note that if we set $a=kn$ and $b=n$
in (\ref{eq:classic}) and (\ref{eq:classicnw}),
then we get the same answers as if we set $s=k$ and $t=1$ in
(\ref{eq:totalone}) and (\ref{eq:nwone}).
(When $t=1$, the map $S\mapsto S^+$
simply adds a northwest corner to every path.)

Our proof of Theorem~\ref{thm:corners} is similar to
Goulden and Serrano's in several ways but differs in one crucial way.
Like Goulden and Serrano, we interpret
(\ref{eq:nwone}) and (\ref{eq:nwtwo}) as counting
{\em all\/} paths of a certain type,
minus the {\em bad\/} paths.
Another similarity is
the idea of breaking the bad path into two halves $\rho$ and~$\sigma$
at the first ``bad point'' so as to
manipulate $\rho$ and~$\sigma$ into something that is easier to count.
The crucial difference is that Goulden and Serrano
{\em rotate\/} $\rho$, whereas we
{\em interchange $\rho$ and~$\sigma$}.\footnote{In
an earlier draft of this paper, we stated,
``Therefore our bijection does not specialize to
Goulden and Serrano's rotation principle nor to
Andr\'e's reflection principle.''
There are many problems with this remark.
First, we learned from Marc Renault~\cite{Re},
Heinrich Niederhausen, and Katherine Humphreys
that Andr\'e did \emph{not} use a reflection argument,
but actually interchanged $\rho$ and~$\sigma$!
We also overstated the difference between rotation and interchange;
if we rotate $\rho$, then rotate $\sigma$,
and then rotate the entire path, then we have
simply interchanged $\rho$ and~$\sigma$.
So when there is a proof by one technique,
there is probably a proof by the other.
Finally, independently and almost simultaneously
with Goulden and Serrano, Loehr~\cite{Loe} used a rotation argument
for a similar lattice-path enumeration problem.}

We also give a second proof of Corollary~\ref{cor:total},
which is based on a well-known argument of Raney~\cite{Ra}
regarding cyclic shifts of integer sequences.

It is frustrating that Theorem~\ref{thm:corners} applies only to
special endpoints.  Can anything be said about other endpoints?
We do not have a satisfactory answer to this question,
but our second main result is a tantalizing hint that
more general theorems lie waiting to be found.
It is best stated in the language of binary strings;
we draw the connection to lattice paths afterwards.

\begin{thm}
\label{thm:binary}
For $n\ge 1$, $s\ge 0$, and $0\le r\le 2n$,
let $a(n,s,r)$ be the number of
binary sequences of length $(s+2)n+1$
such that for all~$j$, the $j$th occurrence of $10$ (if it exists)
appears in positions $(s+2)j+1$ and $(s+2)j+2$ or later,
and such that the total number of
occurrences of $10$ and~$01$ is at most~$r$.  Then
\begin{equation}
\label{eq:binary}
a(n,s,r) = 2\binom{(s+2)n-1}{r} - (s-2) \sum_{i=0}^{r-1}
   \binom{(s+2)n-1 }{i}.
\end{equation}
\end{thm}

Our proof of Theorem~\ref{thm:binary} is again an
application of Raney's argument,
combined with a straightforward induction on~$n$.

To convert Theorem~\ref{thm:binary} into lattice-path language,
let $\beta = (b_1, b_2, \ldots, b_{(s+2)n+1})$ be a binary sequence,
let $b_0 = 0$, and define
$\Delta \beta$ by $(\Delta \beta)_i = |b_i - b_{i-1}|$, for $i\ge1$.
If we convert $\Delta \beta$ into a lattice path by
turning $0$'s into east steps and $1$'s into north steps,
then it is easily checked that the binary sequences in
Theorem~\ref{thm:binary} turn into lattice paths avoiding~$B_s$,
as defined below.

\begin{defn}
For $s\ge 0$, define $B_s$ to be the staircase path that starts at $(0,2)$,
then takes $s+1$ steps east, $2$ steps north, $s$ steps east,
$2$ steps north, $s$ steps east, and so on, always alternating
between $2$ steps north and $s$ steps east except for the first
segment of $s+1$ steps east.
\end{defn}

For example, $B_4$ is the dashed line in the lower picture
in either Figure~\ref{fig:caseone} or Figure~\ref{fig:casetwo} below.
Curiously, we have not been able to generalize Theorem~\ref{thm:binary}
to more general staircase boundaries,
or to refine the count according to northwest or southeast corners.
But for the special case when $s=2k$ is even,
we have a second, purely bijective proof of the following
corollary of Theorem~\ref{thm:binary}.

\begin{cor}
\label{cor:sum}
For all $n\ge1$ and $k \ge 0$,
the number of lattice paths of length $2(k+1)n+1$
that start at $(0,0)$ and that avoid touching or crossing $B_{2k}$
equals the number of lattice paths of length $2(k+1)n+1$ that start 
at $(0,0)$ and that avoid touching or crossing the line $x=ky$ except
at $(0,0)$.  This number has the explicit formula
\begin{equation}
\label{eq:sum}
\binom{2(k+1)n }{2n} - (k-1) \sum_{i=0}^{2n-1}
   \binom{2(k+1)n }{i}.
\end{equation}
\end{cor}

Formula (\ref{eq:sum}) is of course just obtained by summing over
the appropriate instances of (\ref{eq:classic}).
These numbers also appear as A107027 in Sloane's
Online Encyclopedia of Integer Sequences.
This cries out for a combinatorial interpretation of each summand
as counting lattice paths avoiding $B_{2k}$ but with varying endpoints.
Unfortunately, we do not know how to make this idea work.

Note that the case $s=2$ of Theorem~\ref{thm:binary} is particularly simple:

\begin{cor}
\label{cor:fournchoosetwon}
For $n\ge 1$, there are $\binom{4n}{2n}$ binary sequences of length $4n+1$
with the property that for all~$j$, the $j$th occurrence of $10$
appears in positions $4j+1$ and $4j+2$ or later
(if it exists at all).
\end{cor}

We suspect that we have not yet found 
the ``proof from the Book'' of
Corollary~\ref{cor:fournchoosetwon},
and encourage the reader to find it.

The outstanding open question is whether our results generalize further.
We should mention two papers \cite{MN} and~\cite{T}
that consider staircase boundaries similar to~$A_{s,t}$
and that prove results related to
Corollary~\ref{cor:total}.
See also Theorem~8.3 of~\cite{BMN}.
Although our results do not seem to imply or be implied by
these other results, perhaps it would be fruitful to investigate the
precise relationships among them.

We thank
Mihai Ciucu,
Don Coppersmith,
Ira Gessel,
Christian Krattenthaler,
Fred Kochman,
Lee Neuwirth,
and
Doron Zeilberger for useful discussions.

\section{Proofs of Theorem \ref{thm:corners} and Corollary \ref{cor:total}}
\label{sec:proofcorners}

\begin{pf} (of Theorem \ref{thm:corners})

We prove part (1) first.
It will be convenient to first prove formula~(\ref{eq:totalone})
bijectively, and then track corner counts through the bijection.

As we hinted above, we interpret (\ref{eq:totalone}) as counting
the set~$T$ of {\em all\/} paths of a certain type,
minus the set of {\em bad\/} paths.
For $0\le i\le t-1$, let $T_i$ be the set of all lattice paths from
$(1,i)$ to $(sn+1,tn+i)$, and let $T=\bigcup_i T_i$.  Then
\begin{equation}
|T_i| = \binom{sn+tn}{tn} \qquad \mbox{and}\qquad
|T| = t\binom{sn+tn}{tn}.
\end{equation}
We regard $S_1$ as a subset of~$T$ as follows.
Given any $\pi\in S_1$, find the smallest~$i$
such that $(1,i)\in\pi$; such an~$i$ must exist.
Then there exists a unique $\pi'\in T_i$ that agrees exactly
with the remainder of~$\pi$, provided that we append $i$ north steps
to the end of~$\pi$.  Identifying $\pi$ with~$\pi'$ embeds $S_1$ in~$T$.
It remains to show that the number of bad paths---i.e.,
the paths in $T\backslash S_1$---is $s\binom{sn+tn}{tn-1}$.

We partition the set $T\backslash S_1$ into $s$ disjoint sets
$U_1, \ldots, U_s$ as follows.
By definition, every path in $T\backslash S_1$ must hit a bad point,
i.e., a point on the boundary $A_{s,t}$.
For $1\le j \le s$,
we let $U_j$ be the set of all paths in $T\backslash S_1$
whose first bad point has an $x$-coordinate that is
congruent to~$j$ modulo~$s$.
To prove formula~(\ref{eq:totalone}), it suffices to show that
$|U_j| = \binom{sn+tn }{tn-1}$, independent of~$j$.

Fix any $j$.  Given $\pi \in U_j$, observe that the step that
terminates in the first bad point of~$\pi$ must be a north step.
Let $\rho$ be the portion of~$\pi$ prior to this fatal north step,
and let $\sigma$ be the portion of~$\pi$ after the bad point.
Thus $\pi = (\rho, \mbox{north}, \sigma)$.
Now comes the crucial part of the proof,
where we interchange $\rho$ and~$\sigma$.
More precisely, let $\pi'$ be the lattice path that
starts at $(j-1,t)$ and takes steps $(\sigma, \mbox{east}, \rho)$.
See Figure~\ref{fig:interchange} for an example.

\begin{figure} [ht]
\centering
\includegraphics{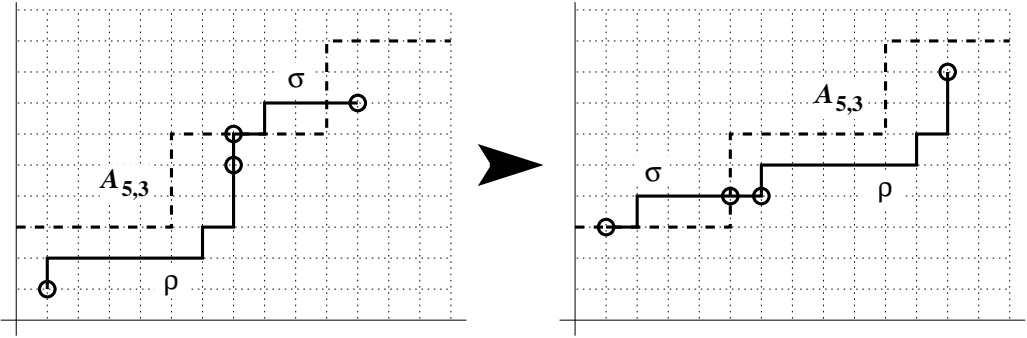}
\caption{Example of $\pi\mapsto\pi'$
with $s=5$, $t=3$, $n=2$, and $\pi \in T_1 \cap U_2$}
\label{fig:interchange}
\end{figure}

We claim that $\pi\mapsto \pi'$ bijects $U_j$ onto
the set~$U_j'$ of all paths from
$(j-1,t)$ to $(sn+j, tn+t-1)$.
First note that since $\pi$ and $\pi'$ have
the same total number of north steps and
the same total number of east steps except that
one north step of~$\pi$ has been changed into an east step of~$\pi'$,
it follows that $\pi'$ does in fact terminate at $(sn+j, tn+t-1)$.
Now, given any path~$\pi'\in U_j'$,
let $\sigma$ be the initial segment of~$\pi'$ up to the last
point of~$\pi'$ that lies on the boundary $A_{s,t}$.
The next step after that must be an east step;
let $\rho$ be the remainder of~$\pi'$ after that.
It is straightforward to check that this allows us to
construct a unique preimage~$\pi$ of~$\pi'$.
This proves formula~(\ref{eq:totalone}).

Now we prove the stronger formula~(\ref{eq:nwone}), for northwest corners.
The initial north step of each lattice path in~$S_1^+$
forces there to be a northwest corner with $x$-coordinate zero,
whereas the final north step does not affect the northwest corner count.
Therefore if we embed $S_1$ in~$T$ as above,
we really want to count lattice paths with $c-1$ northwest corners
(rather than $c$ northwest corners).
There are $t \binom{sn}{c-1}\binom{tn}{c-1}$
paths in~$T$ with $c-1$ northwest corners,
because we can pick the $x$-coordinates and $y$-coordinates
of the corners independently.
It therefore suffices to show that for all~$j$,
there are $\binom{sn-1 }{c-2} \binom{tn+1}{c}$ paths in~$U_j$
with $c-1$ northwest corners.

If $\alpha$ is a binary string,
let $|\alpha|$ denote its length,
and let $w(\alpha)$ denote its {\em weight,}
i.e., the number of $1$'s in~$\alpha$.
Let $U_j''$ be the set
of ordered pairs $(\alpha,\beta)$ of binary strings
such that $|\alpha| = sn-1$, $|\beta| = tn+1$,
and $w(\beta) = w(\alpha) + 2$.
It suffices to describe a bijection from $U_j'$ to $U_j''$
such that the composite map $\pi \mapsto \pi' \mapsto (\alpha,\beta)$
sends paths with $c-1$ northwest corners to pairs $(\alpha,\beta)$
with $w(\beta)=c$.

Before describing this bijection, we make two observations.
Let $\pi$, $\rho$, $\sigma$, and~$\pi'$ be as above.
The first observation is that,
because of the position of the endpoint of~$\pi$
relative to the boundary $A_{s,t}$,
$\sigma$ always has at least one east step.
The second observation is that
we lose a northwest corner when passing from $\pi$ to~$\pi'$
iff $\sigma$ starts with an east step,
and we gain a northwest corner as we pass from $\pi$ to~$\pi'$
iff $\sigma$ ends with a north step.
(Note that we can both gain a corner and lose a corner,
leaving the total corner count unchanged.)
So to track corners properly,
we must watch the first and last steps of~$\sigma$.

Now for the bijection.
Given $\pi'\in U_j'$, construct $\alpha$ by first writing down
a binary string of length $sn$ whose $i$th digit ($1\le i\le sn$) is~$1$
iff $j+i-1$ is the $x$-coordinate of a northwest corner of~$\pi'$,
and then deleting the digit corresponding to the point where $\pi'$
intersects $A_{s,t}$ for the last time.
This digit must exist, because $\sigma$ has at least one east step.
For example, in Figure~\ref{fig:interchange},
we first write down $1000100001$, and then delete the $4$th digit
to obtain $\alpha=100100001$.

The first $tn-1$ digits of $\beta$ are obtained by writing down
the binary string of 
length $tn-1$ whose $i$th digit ($1\le i\le tn-1$) is~$1$
iff $t+i$ is the $y$-coordinate of a northwest corner of~$\pi'$.
The next digit of~$\beta$ is~$1$ iff
$\sigma$ does not start with a north step,
and the last digit of~$\beta$ is
the complement of the deleted digit of~$\alpha$.
For example, in Figure~\ref{fig:interchange},
$\beta=1110011$.

To see that $w(\beta) = w(\alpha)+2$,
first pair off the $1$'s in $\alpha$ and~$\beta$
arising from northwest corners that they both ``see,''
and then note that $\beta$ will have two extra $1$'s
corresponding to the 
columns in which the first and last vertices of~$\sigma$ appear:
Either $\beta$ sees a northwest corner in that column
(and $\alpha$ of course does not see it),
or there is no such corner, in which case
the appropriate trailing bit of~$\beta$ will be set.
Either way, $w(\beta)=w(\alpha)+2$.

Similarly, as we pass from $\pi$ to $\pi'$ to~$\beta$,
a corner that is lost from $\pi$ to~$\pi'$ is ``caught'' by
the penultimate bit of~$\beta$,
and $\beta$ will gain an extra~$1$ either by catching
a gained corner or, if no corner is gained, by setting its last bit.
Thus $w(\beta)$ is one more than the number of northwest corners of~$\pi$.
Equivalently, $w(\alpha)$ is one {\em less\/} than
the number of northwest corners of~$\pi$.

It remains to show that $\pi'\mapsto(\alpha,\beta)$ is a bijection.
Since $|U_j'| = |U_j''|$, it suffices to show
that $\pi'$ can be reconstructed from its image $(\alpha,\beta)$.
To reconstruct~$\pi'$ it suffices to reconstruct the northwest corners.
The penultimate digit of~$\beta$ is~$0$
iff $\pi'$ has a northwest corner with $x$-coordinate $j-1$,
so we need only reconstruct the deleted digit of~$\alpha$.
The {\em value\/} of the deleted digit is the complement of
the last digit of~$\beta$,
so we need only reconstruct its {\em position\/}.
To do this, take $(\alpha,\beta)$ and begin constructing $\pi'$ from
the end backwards without regard to the deleted digit.
At some point, the partially reconstructed path will touch or cross
the boundary $A_{s,t}$.
It is easy to check that the first such contact point with $A_{s,t}$
yields the position of the deleted digit of~$\alpha$.

This completes the proof of part~(1).
The proof of part~(2) is very similar,
so we focus only on the details that differ.
For $0\le i\le t-1$, let $T_i$ be the set of all lattice paths from
$(1,i)$ to $(sn,tn+i-1)$, and let $T=\bigcup_i T_i$.
Then $T$ is our set of {\em all\/} paths.
Note that $S_2$ is already
naturally a subset of~$T$---in fact, $S_2 \subset T_0$---so
we do not have to embed $S_2$ in~$T$.
The definition of the sets $U_j$ is exactly analogous.
However, $\pi'$ now starts at $(j,t+1)$ rather than at $(j-1,t)$,
and ends at $(sn+j, tn+t-1)$.
The proof of~(\ref{eq:totaltwo}) now goes through as before.

To do the corner count,
we need to define the map $\pi'\mapsto (\alpha, \beta)$ in
the case that $\sigma$ is empty or vertical, i.e., has no east steps.
In this case, we always delete the {\em first\/} digit of~$\alpha$.
The definition of~$\beta$ is the same as before.
The arguments that $w(\beta) = w(\alpha) + 2$ and that
$\pi'\mapsto (\alpha, \beta)$ is a bijection still work.

However, $w(\alpha)$ is no longer always one less than the
number of northwest corners of~$\pi$.
Let $V$ denote the set of paths in~$\bigcup_j U_j$
for which $\sigma$ is vertical or empty
and $\rho$ starts with a horizontal step.
Then it is straightforward to check that for $\pi\in V$,
$w(\alpha)$ is {\em equal\/} to the number of northwest corners of~$\pi$.
So if we let $X^c$ denote the members of~$X$ with $c$ northwest corners,
then pulling back $U_j''$ to $U_j$ shows that
formula~(\ref{eq:nwtwo}) is the cardinality of the set
\begin{equation}
\label{eq:parttwoa}
\biggl(T^{c-1} \backslash \bigcup_{j=1}^s U_j^{c-1} \biggr) \cup
  (V^{c-1} \backslash V^{c-2})
= S_2^{c-1} \cup
  (V^{c-1} \backslash V^{c-2}).
\end{equation}
On the other hand, if we let $N_2$ denote the subset of~$S_2$
consisting of paths that start with a north step,
and observe that prepending a north step to $\pi\in S_2$
adds a northwest corner to~$\pi$ iff $\pi$ starts with an east step,
then we see that $(S_2^+)^c$ is equinumerous with
\begin{equation}
\label{eq:parttwob}
S_2^{c-1} \cup
  (N_2^{c} \backslash N_2^{c-1}).
\end{equation}
Thus to show that (\ref{eq:parttwoa}) and (\ref{eq:parttwob})
are equinumerous,
it suffices to show that $N_2^c = V^{c-1}$ for any~$c$.
But this bijection is easily described:
Given a path in~$N_2$, simply move all the initial north steps
to the end; this creates a path in~$V$ with one fewer northwest
corner.  This completes the proof.
\end{pf}

\begin{pf} (of Corollary~\ref{cor:total})

Of course this follows from Theorem~\ref{thm:corners},
but we have another proof.
The formula in equation~(\ref{eq:totalone})
can be rewritten as $\frac{1}{n}\binom{sn+tn}{sn+1}$.
Consider the set~$S_1'$ of all paths starting at the origin
that end with a north step and
that have a total of $sn+1$ east steps
and a total of $tn$ north steps.
Clearly $|S_1'| = \binom{sn+tn}{sn+1}$ and $S_1 \subset S_1'$.
We need to show that $S_1$ comprises precisely $1/n$ of the paths in~$S_1'$.

Decompose any path $\pi\in S_1'$ into $n$ consecutive subpaths
$\pi_1, \pi_2, \ldots, \pi_n$,
where each $\pi_j$ contains exactly $t$ north steps
and ends in a north step.
Our desired result follows immediately from the following key claim:
For {\em any\/} $\pi\in S_1'$,
there is exactly one ``cyclic shift'' of $\pi$ that lies in~$S_1$,
where by a cyclic shift we mean one of the $n$ paths
of the form
$$\pi_j, \pi_{j+1}, \ldots, \pi_n, \pi_1, \pi_2, \ldots, \pi_{j-1}$$
obtained from~$\pi$ by concatenating the subpaths in a cyclically
permuted order.

To see the key claim,
one first readily verifies that $\pi\in S_1$ iff for all $i>0$,
the total length of the first $i$ subpaths
$\pi_1, \ldots, \pi_i$ is at least $(s+t)i + 1$.
Now we apply an argument patterned after a classic proof of Raney~\cite{Ra}.
For all $j\ge 1$, let $\ell_j$ be the length of~$\pi_{(j\bmod n)}$.
Consider the graph in the $xy$ plane with straight-line segments
between vertices $P_j$ and $P_{j+1}$, where
$$P_j = \biggl(j, \sum_{i=1}^j \ell_i\biggr).$$
The ``average'' slope of this graph is $(sn+tn+1)/n = s+t+\frac{1}{n}$.
The line of the form $y=(s+t+\frac{1}{n})x + C$
that is ``tangent'' to this graph from below
intersects the graph exactly once every $n$ points,
because the graph has period~$n$ and
the coefficient of~$x$ is an integer plus $1/n$.
The points of intersection have the form $P_j$, $P_{j+n}$, $P_{j+2n}$, etc.,
and the value of~$j$ here yields the unique cyclic shift
having the desired property.  This proves the claim.

The proof of equation~(\ref{eq:totaltwo}) is similar.
Define $S_2'$ to be the set of all paths from the origin
that have a total of $sn-1$ east steps and a total of $tn-1$ north steps.
Decompose any $\pi\in S_2'$ as follows:
$$\pi = \pi_1, \mbox{north}, \pi_2, \mbox{north}, \ldots, \mbox{north},
  \pi_n$$
where each $\pi_j$ has $t-1$ north steps.
Then $\pi\in S_2$ iff for all $0\le i<n$ we have
$$|\pi_1| + |\pi_2| + \cdots + |\pi_i| \ge i(s+t-1).$$
Exactly one ``cyclic shift'' of~$\pi$ has the equivalent property that
for all~$i$,
$$|\pi_1| + |\pi_2| + \cdots + |\pi_i| \ge (i/n)(sn+tn-n-1).$$
Thus there are $\frac{1}{n} \binom{sn+tn-2}{tn-1}$ paths in~$S_2$,
which is equivalent to equation~(\ref{eq:totaltwo}).
\end{pf}

\section{Proofs of Theorem \ref{thm:binary} and Corollary \ref{cor:sum}}
\label{sec:proofmain}

\begin{pf} (of Theorem \ref{thm:binary})

It is easily verified that $a(1,s,0)=2$ and
$a(1,s,1)=a(1,s,2)=s+4$. If $n \ge 2$ and $r \le 2n-2$, then we
claim that the following recursion holds:
\begin{equation}
\label{eq:recursion}
a(n,s,r) = \sum_{d=0}^{s+2} \binom{s+2}{d} a(n-1,s,r-d).
\end{equation}
The reason is that an admissible binary string of order $n-1$ can
be extended by any sequence of $s+2$ bits without danger of
causing inadmissibility, provided that the resulting string
changes from 1 to 0 or vice versa at most $2n-2$ times. The
parameter $d$ counts the number of changes introduced by the
last $s+2$ bits, and the binomial coefficient counts the number
of ways to position the $d$ changes.

By Vandermonde convolution,
the recurrence (\ref{eq:recursion})
almost gives us a proof by induction on~$n$,
except that we need to handle the cases $r=2n-1$ and $r=2n$.
Note that no string of order $n$ can have more than $2n-1$
changes, and that equation~(\ref{eq:binary})
takes the same value for $r=2n-1$
and $r=2n$. So to complete the proof of Theorem~\ref{thm:binary},
it is enough to show that
\[
a(n,s,2n-1) - a(n,s,2n-2) =
2 \binom{(s+2)n-1}{2n-1} - s \binom{(s+2)n-1}{2n-2} ,
\]
which can be rewritten as
$\frac{1}{n} \binom{(s+2)n}{2n-1}$.
The left-hand side counts the admissible strings with exactly
$2n-1$ changes, and we use the proof technique of Raney as before.
Any such string $\sigma$ must start with 0; we decompose it into
substrings $\sigma_1,\sigma_2,\ldots,\sigma_n$, where each
$\sigma_j$ consists of $a_j$ zeroes followed by $b_j$ ones,
and $a_j,b_j>0$. The condition for admissibility can now be
expressed as
\[
|\sigma_1| + |\sigma_2| + \cdots + |\sigma_i| \geq i(s+2)+1
\]
for all $0 \leq i < n$. Exactly one cyclic shift of~$\sigma$
% $\sigma_j,\sigma_{j+1},\ldots,\sigma_n,%
% \sigma_1,\ldots,\sigma_{j-1}$
has the equivalent property that for all $i$,
\[
|\sigma_1| + |\sigma_2| + \cdots + |\sigma_i| \geq (i/n)((s+2)n+1).
\]
Thus the number of admissible strings with exactly $2n-1$ changes
is equal to $1/n$ times the number of ways to partition
$(s+2)n+1$ into $2n$ positive integers, corresponding to the
numbers $a_j$, $b_j$. This is well known to be
$\binom{(s+2)n}{2n-1}$, and this completes the proof.
\end{pf}

\begin{pf} (of Corollary \ref{cor:sum})

We can deduce this easily from Theorem~\ref{thm:binary} 
just by showing that equation~(\ref{eq:binary})
reduces to equation~(\ref{eq:sum}) when $r=2n$ and $s=2k$.  We have
\[
a(n,2k,2n) = 2 \binom{2(k+1)n-1}{2n}
- 2(k-1) \sum_{i=0}^{2n-1} \binom{2(k+1)n-1}{i} ;
\]
breaking up the sum, the right-hand side becomes
\begin{align*}
& 2 \binom{2(k+1)n-1}{2n}
- (k-1) \binom{2(k+1)n-1}{2n-1}
\\ & \quad
- (k-1) \sum_{i=0}^{2n-1}
\left[ \binom{2(k+1)n-1}{i} + \binom{2(k+1)n-1}{i-1} \right] ,
\end{align*}
or
\[
\left[ 2 \frac{2kn}{2(k+1)n} -
(k-1) \frac{2n}{2(k+1)n} \right]
\binom{2(k+1)n}{2n}
- (k-1) \sum_{i=0}^{2n-1} \binom{2(k+1)n}{i},
\]
which then collapses to formula~(\ref{eq:sum}).

However, we also give a direct bijective proof.
If $k=0$ then formula~(\ref{eq:sum}) simplifies to~$4^n$,
the boundary conditions are nearly vacuous,
and the result is easy to prove.
So fix $k\ge 1$ and $n\ge 1$.

Our bijection is actually between two sets of lattice paths that are slightly
different from those mentioned in the corollary.

Let $\fkn$ be the set of lattice paths of length $2(k+1)n$
(note the shorter length)
that start at $(0,0)$ and avoid the line $x=ky$ except at $(0,0)$.

Let $\gkn$ be the set of lattice paths of length $2(k+1)n+1$
that start at $(0,0)$ and avoid~$B_{2k}$, and that touch
the line $x=ky+1$ at least once for $y>0$.

To see that a bijection between $\fkn$ and $\gkn$ implies the corollary,
we make two observations.
First, because $2(k+1)n$ is a multiple of~$k+1$,
every lattice path of length $2(k+1)n$
that avoids $x=ky$ can be extended by
either an east step or a north step without hitting the line $x=ky$;
therefore $|\fkn|$ is exactly half the number of lattice paths
of length $2(k+1)n+1$ that avoid $x=ky$.
Second, the paths excluded by the final condition on~$\gkn$
are precisely those that avoid the line $x=ky+1$ after $(1,0)$,
and therefore are in bijection with $\fkn$---simply prepend
an east step to each path in~$\fkn$.

The rest of the proof is devoted to describing a bijection
$\varphi: \fkn \to \gkn$.

We define a procedure called {\em trisection\/} that we need 
in our construction of~$\varphi$.
Define the {\em potential\/} of a point $(x,y)$ to be $x-ky$.
Let $P$ be a lattice path, not necessarily starting at $(0,0)$,
but with the property that the {\rm potential difference\/}
of~$P$---i.e., the potential of the last point of~$P$ minus
the potential of the first point of~$P$---is at least $k+1$
(and hence in particular comprises at least $k+1$ steps).
To trisect~$P$, first look at the last $k$ steps of~$P$.
If all of these steps are east steps, then the trisection procedure fails.
Otherwise, let $b$ be the segment of~$P$ consisting of the
last north step of~$P$ along with all the east steps after that.
Let $l$ be the length of~$b$.
Find the {\em last\/} lattice point $p\in P$
such that the initial segment~$a$ of~$P$
comprising everything up to~$p$
has potential difference exactly $k+1-l$.
Such a point $p$ must exist
(since increases in potential can occur
only one unit at a time) and must occur prior to~$b$.
Let $P'$ be the segment of~$P$ between $a$ and~$b$.
The decomposition $P = (a, P', b)$ is the trisection of~$P$.
Note that the potential difference of~$b$ is $l-1-k$
and so the combined potential difference of $a$ and~$b$ is zero;
thus the potential difference of~$P'$ is the same as that of~$P$.

We are now ready to describe $\varphi$.
Given $P \in \fkn$,
the construction of $\varphi(P)$ has two phases.
In Phase~1, we decompose $P$ into segments;
in Phase~2, we build $\varphi(P)$ using the
segments constructed in Phase~1.

The paths in~$\fkn$ with the smallest potential difference
are those that terminate closest to the line~$x=ky$;
these are readily checked to have potential difference at least $k+1$.
We begin Phase~1 by trying to trisect $P$ into $(a_1, P', b_1)$.
If this fails, we proceed to Phase~2.
Otherwise, if the {\em height\/} of~$a_1$
(i.e., the $y$-coordinate of the last point of~$a_1$
minus the $y$-coordinate of the first point of~$a_1$)
is even, then we proceed to Phase~2.
Otherwise, we try to trisect $P'$ into $(a_2, P'', b_2)$,
proceeding to Phase~2 if the trisection fails or if $a_2$ has even height.
If we still do not reach Phase~2,
then we try to trisect~$P''$, and so~on.

Each successful trisection preserves the potential difference
of the middle section while shrinking its length,
so we must eventually reach Phase~2, with a decomposition
$$P = (a_1, a_2, \ldots, a_{m-1}, a_m, Q, b_m, b_{m-1},
  \ldots, b_2, b_1),$$
for some~$m$, where $Q$ denotes whatever remains in the middle.
If we reach Phase~2 because the height of~$a_m$ is even, then we set
$$ \varphi(P) = (\mbox{east}, a_m, b_m, a_1, b_1, a_2, b_2, a_3, b_3,
   \ldots, a_{m-1}, b_{m-1}, Q), $$
where the ``east'' means that we begin $\varphi(P)$ with an east step.
For an example with $k=2$ and $n=7$, see Figure~\ref{fig:caseone}.

\begin{figure} [ht]
\centering
\includegraphics{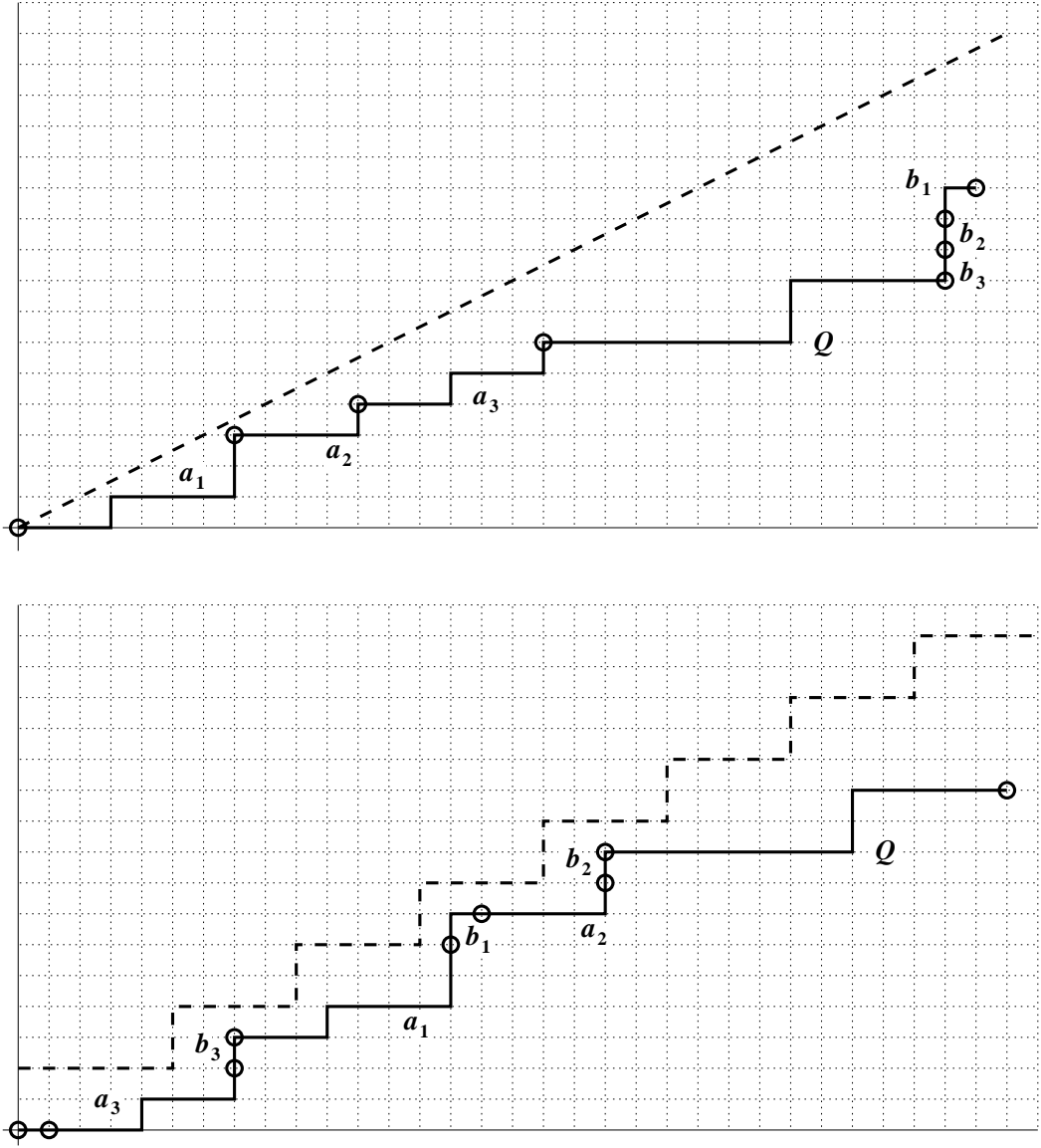}
\caption{Example of $P \mapsto \varphi(P)$ when $a_m$ has even height}
\label{fig:caseone}
\end{figure}

The other way to reach Phase~2 is
for the last $k$ steps of~$Q$ to all be east steps.
Decompose $Q = (Q', b_{m+1})$ where $b_{m+1}$ comprises those final
$k$ east steps.  Then set
$$ \varphi(P) = (\mbox{north}, \mbox{east}, b_{m+1}, a_1, b_1, a_2, b_2,
   \ldots, a_m, b_m, Q').$$
For an example, again with $k=2$ and $n=7$, see Figure~\ref{fig:casetwo}.

\begin{figure} [ht]
\centering
\includegraphics{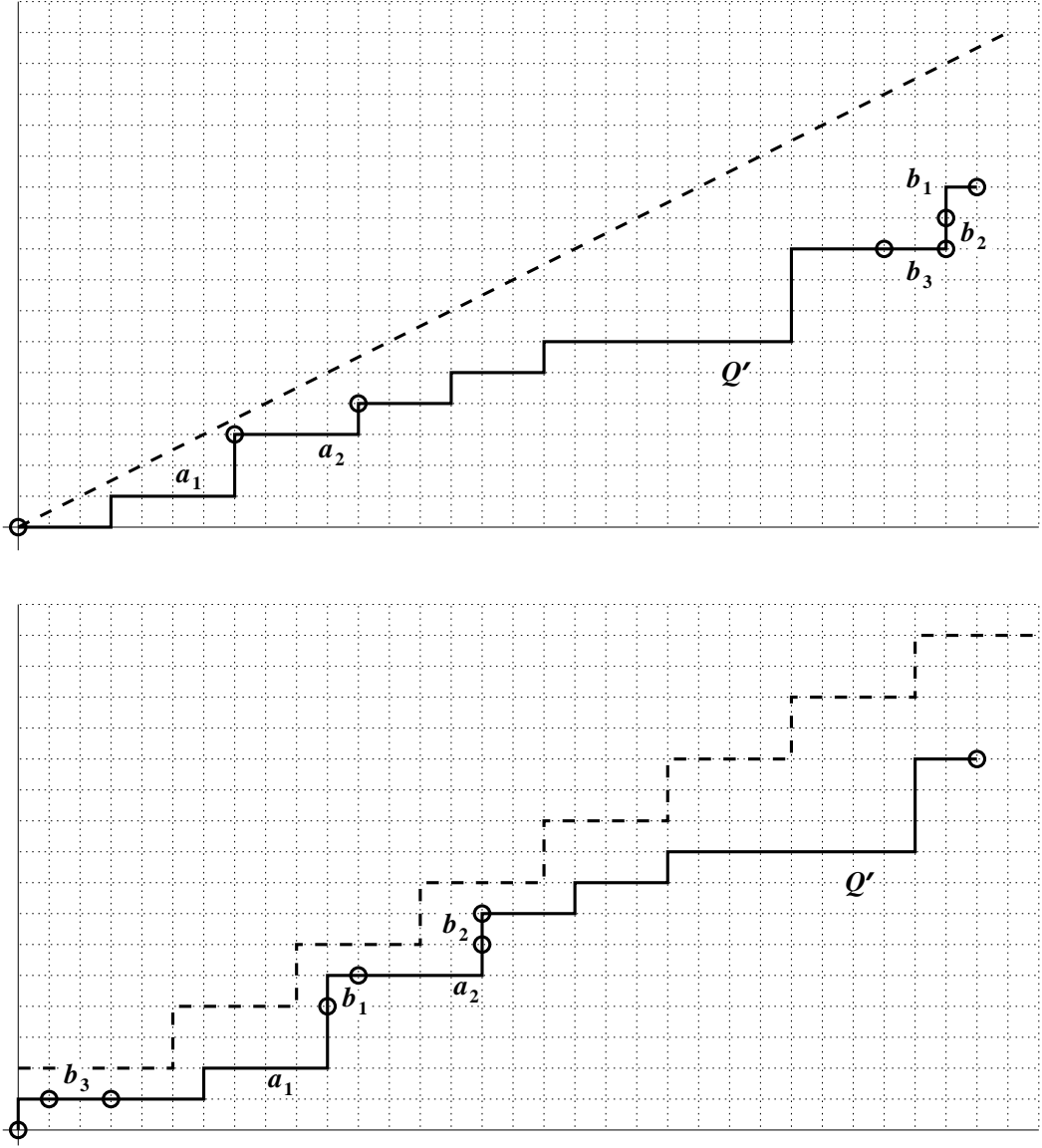}
\caption{Example of $P \mapsto \varphi(P)$ when $b_{m+1}$ comprises $k$
east steps}
\label{fig:casetwo}
\end{figure}

One must check that $\varphi(P) \in \gkn$.
Define a \emph{waypoint} to be a point
on the line $x=ky+1$ that is midway between
two consecutive southeast corners of~$B_{2k}$
(i.e., it has the form $(2ik+k+1, 2i+1)$ for some~$i$).
The claim that $\varphi(P) \in \gkn$ can be broken down into subclaims:
(1)~the part of $\varphi(P)$ preceding~$a_1$
avoids $B_{2k}$ and takes us to a waypoint;
(2)~each pair $(a_1,b_1)$, $(a_2,b_2)$, etc.,
resulting from successful trisections
takes us from one waypoint to another and avoids~$B_{2k}$;
(3)~the last part $Q$ or~$Q'$ starts at a waypoint
and avoids~$B_{2k}$.
Checking these subclaims is easier to do oneself
than to write out in detail,
so we will just indicate the key points.
The claims about waypoints follow because
an $(a_i,b_i)$ pair has potential difference zero,
and therefore if it starts on $x=ky$ then it ends on $x=ky$.
If $a_i$ has odd height then $(a_i,b_i)$ ends on a waypoint
if it starts on a waypoint.
At most one~$a_i$, namely~$a_m$, has even height,
and then $(a_m,b_m)$ takes us from $(1,0)$ to a waypoint.
The other tricky claim is that
$(a_i,b_i)$ avoids~$B_{2k}$,
but this follows because by construction,
the potential of $a_i$
never drops below its initial potential
so that it even stays below the line $x=ky+1$;
also, if we trace $b_i$ backwards from its terminal waypoint,
it takes at most $k-1$ horizontal steps and therefore avoids
hitting~$B_{2k}$.

To invert~$\varphi$, suppose we are given $P\in\gkn$.
Whether Figure~\ref{fig:caseone} or Figure~\ref{fig:casetwo}
applies depends on whether the first step of~$P$ is north or east.
Mark all the waypoints of~$P$;
there must be at least one, since $P\in\gkn$.
By backing up from a waypoint until we find a north step,
we can construct the~$b_i$,
and therefore also the~$a_i$ and~$Q$.
Hence $\varphi$ is easily reversed.
We leave the straightforward verification of the
details to the reader.
\end{pf}

% For example, we might ask whether there is a simple formula
% for the number of lattice paths from $(0,0)$ to $(x,y)$
% that avoid~$B_k$.
% For $k=1$, we conjecture that this quantity is
% $$\sum_{i=0}^y a_i \binom{x+y-i}{y-i},$$
% where
% $$\sum_{i\ge 0} 2 a_i z^i = (1-2z)(1+4z-4z^2) + (1-2z)^2 \sqrt{1+4z-4z^2}.$$
% % If $x=2n+3$ and $y=2n+1$ then a simpler expression appears possible,
% % namely $4^{n+1} C_n$ where $C_n$ is the $n$th Catalan number.
% For larger~$k$, however, we do not know what these numbers are.

% \section{Acknowledgments}

% \appendix

% \section{}
% \label{}

\end{document}